\documentclass[12pt, a4paper]{amsart}
\usepackage{amscd, amssymb,amsmath}
\allowdisplaybreaks[2]
\begin{document}

\def\supp{\operatorname{supp}}
\def\Ind{\operatorname{Ind}}
\def\Aut{\operatorname{Aut}}
\def\Ad{\operatorname{Ad}}
\def\range{\operatorname{range}}
\def\sp{\operatorname{sp}}
\def\clsp{\overline{\operatorname{sp}}}
\def\Prim{\operatorname{Prim}}
\def\dashind{\operatorname{\!-Ind}}
\def\id{\operatorname{id}}
\def\rt{\operatorname{rt}}
\def\lt{\operatorname{lt}}

\def\H{\mathcal{H}}
\def\L{\mathcal{L}}
\def\K{\mathcal{K}}
\def\D{\mathcal{D}}
\def\C{\mathbb{C}}
\def\T{\mathbb{T}}
\def\Z{\mathbb{Z}}
\def\Q{\mathcal{Q}}

\def\m{\mathfrak{m}}
\def\n{\mathfrak{n}}
\def\p{\mathfrak{p}}

\newtheorem{thm}{Theorem}  [section]
\newtheorem{cor}[thm]{Corollary}
\newtheorem{lemma}[thm]{Lemma}
\newtheorem{prop}[thm]{Proposition}
\newtheorem{thm1}{Theorem}

\theoremstyle{definition}
\newtheorem{definition}[thm]{Definition}
\newtheorem{remark}[thm]{Remark}
\newtheorem{example}[thm]{Example}
\newtheorem{remarks}[thm]{Remarks}
\newtheorem{claim}[thm]{Claim}
\newtheorem{problem}[thm]{Problem}

\numberwithin{equation}{section}

\title[Mansfield's imprimitivity theorem]{\boldmath Mansfield's imprimitivity 
theorem for arbitrary closed subgroups}

\author[an Huef]{Astrid an Huef}
\address{School of Mathematics\\
The University of New South Wales\\
Sydney, NSW 2052\\
Australia}
\email{astrid@maths.unsw.edu.au}

\author[Raeburn]{Iain Raeburn}
\address{School of Mathematical and Physical Sciences\\
University of Newcastle\\
NSW 2308\\ Australia}
\email{iain@frey.newcastle.edu.au}

\thanks{This research was supported by grants from the Australian Research Council, the
University of New South Wales and the University of Newcastle.}


\date{May 7, 2002}

\begin{abstract}
Let $\delta$ be a nondegenerate coaction of $G$ on a $C^*$-algebra $B$, and let $H$
be a closed subgroup of $G$.  The dual action $\hat\delta:H\to\Aut(B\times_\delta
G)$ is proper and saturated in the sense of  Rieffel, and the generalised fixed-point
algebra is the crossed product of
$B$ by the homogeneous space $G/H$. The resulting Morita equivalence is a version
of Mansfield's imprimitivity theorem which requires neither amenability nor normality
of~$H$.
\end{abstract}

\maketitle

\section{Introduction}

Nonabelian duality tells us how to recover information about a $C^*$-dynamical system
$(A,G,\alpha)$ from its crossed product $A\times_\alpha G$. The dual systems consist of
a coaction $\delta:B\to M(B\otimes C^*_r(G))$ of a locally compact group $G$ on a
$C^*$-algebra $B$; one recovers $(A,G,\alpha)$ from $A\times_\alpha G$ by taking the
crossed product $(A\times_\alpha G)\times_{\hat\alpha}G$ by the dual 
coaction $\hat\alpha$, and invoking the duality theorem of Imai and Takai. The 
theory of crossed products by coactions runs parallel to the theory of ordinary 
crossed products by actions, and there are now analogues for coactions of many 
of the important theorems about ordinary crossed products. In particular, 
Mansfield has proved an imprimitivity theorem for crossed products by coactions 
analogous to that of Rieffel and Green for ordinary crossed products \cite{man}.

Let $\delta:B\to M(B\otimes C^*_r(G))$ be a coaction of a locally compact group $G$, and
let
$N$ be a closed subgroup of $G$ which is both normal and amenable. The coaction $\delta$
restricts to a coaction $\delta|$ of $G/N$, and Mansfield's imprimitivity theorem says
that
$B\times_{\delta|}(G/N)$ is Morita equivalent to the crossed product $(B\times_\delta
G)\times_{\hat\delta} N$ by the restriction of the dual action. This theorem 
is widely regarded as the deepest result in the subject, and Mansfield's 
analysis is subtle and complicated. 

The amenability hypothesis in Mansfield's theorem ensures that the restriction $\delta|$
is well-defined. If one considers instead full coactions $\epsilon:B\to M(B\otimes
C^*(G))$, then the restriction $\epsilon|$ makes sense for nonamenable $N$, and
Kaliszewski and Quigg showed that Mansfield's theorem extends to many, but not all, full
coactions
\cite{KQ}. Recently there have been indications that one might be able to extend
Mansfield's theorem to non-normal subgroups \cite{EKR, ahrw-symimp}, and indeed such an
extension has been proved for discrete groups \cite{EQ}. 

Here we  prove a version of Mansfield's imprimitivity theorem which requires
neither amenability nor normality of the subgroup. Our result contains interesting new
information even for amenable non-normal subgroups (Corollary~\ref{cor-amenable}) and for
normal nonamenable subgroups (Corollary~\ref{cor-normal}). We use several powerful 
tools: Rieffel's theory of proper actions \cite{rie-pr, rie-pr2}, the dense 
subalgebras constructed by Mansfield in \cite[\S3]{man}, an averaging technique 
developed by Quigg \cite{quigg, QR}, and the crossed products by homogeneous 
spaces introduced in \cite{EKR}. 

We begin with a short section in which we establish conventions
and notation, and then prove our main Theorem~\ref{mainthm} in 
Section~\ref{guts}.  Our proof is quite different from Mansfield's, and is 
substantially shorter because we are taking advantage of established technology. 
We prove that the dual action $\hat\delta$ of a closed subgroup $H$ on 
$B\times_\delta G$ is proper and saturated with respect to Mansfield's 
subalgebra $\D$, and then identify the associated fixed-point algebra with the 
crossed product of $B$ by the homogeneous space $G/H$. In Section~\ref{normal} 
we look at what happens when $H=N$ is normal.  We use results of Quigg 
\cite{quigg-fr} to define the restriction of a reduced coaction $\delta$ of $G$ 
to a reduced coaction $\delta|$ of $G/N$, identify the crossed product by 
$\delta|$ with the crossed product by the homogeneous space 
(Proposition~\ref{prop-allOK}), and deduce a very satisfactory version of 
Mansfield's theorem for $B\times_{\delta|}(G/N)$. We then reconcile our results 
with those of Kaliszewski and Quigg for full coactions.

\section{Preliminaries}\label{prelim}

Let $G$ be a locally compact group
and $\lambda^G=\lambda$  the left-regular representation of $G$ on $L^2(G)$.
We  view the reduced group $C^*$-algebra $C^*_r(G)$ as the subalgebra
$\lambda(C^*(G))$ of $B(L^2(G))$.
We use minimal tensor products throughout. We routinely extend nondegenerate
homomorphisms to multiplier algebras without changing notation.

In Section~\ref{guts} we use the
conventions of \cite{lprs}. 
Thus $\delta_G$ denotes the
reduced comultiplication
$\delta_G:C^*_r(G)\to M(C^*_r(G)\otimes C^*_r(G))$ characterised by
$\delta_G(\lambda_s)=\lambda_s\otimes\lambda_s$, and  a
\emph{reduced coaction}
$\delta$ of $G$ on a $C^*$-algebra $B$ is an injective nondegenerate
homomorphism $\delta:B\to M(B\otimes C_r^*(G))$ such that
$
(\delta\otimes\iota)\circ\delta=(\iota\otimes\delta_G)\circ\delta 
$ and 
$\delta(b)(1\otimes z)\in B\otimes C_r^*(G)
$
for all $b\in B$ and  $z\in C^*_r(G)$.   
We write $A_c(G)$ for the set of functions in
the Fourier algebra $A(G)$ with compact support.
For $u\in A(G)$ and $b\in B$ we write $\delta_u(b):=(\iota\otimes
u)(\delta(b))$;  the coaction is  \emph{nondegenerate} if
$\overline{\delta_{A(G)}(B)}=B$.
Our main Theorem~\ref{mainthm}, like Mansfield's, is about
reduced coactions. In Section~\ref{normal} we  need to use
\emph{full coactions}, and our conventions there are also the standard ones.

Let $\pi$ be a faithful nondegenerate representation of $B$ on $\H_\pi$ and
$M$ the representation of $C_0(G)$ by multiplication operators on $L^2(G)$.
The closed span
  \[
B\times_{\delta} G:=\clsp\{((\pi\otimes\iota)\circ\delta(b))(1\otimes
M(f)):b\in B, f\in C_0(G)\}.
\]
is a $C^*$-subalgebra of $B(\H_\pi\otimes L^2(G))$, which is up to isomorphism
independent of the representation $\pi$, and which is called the 
\emph{crossed product}
of $B$ by $\delta$. Conjugation by the right-regular representation gives a
\emph{dual action} $\hat\delta$ of $G$ on
$B\times_\delta G$, which is characterised by
\[
\hat\delta_s\big( ((\pi\otimes\iota)(\delta(b)))(1\otimes M(f)) \big)
=((\pi\otimes\iota)(\delta(b)))(1\otimes M(\rt_s(f)))
\]
where $\rt_s(f)(t)=f(ts)$.

Let $H$ be a closed subgroup of $G$.  As in \cite[Definition~2.1]{EKR}, the \emph{reduced
crossed product}
$B\times_{\delta,r} (G/H)$ of $B$ by the
``coaction'' of the homogeneous space $G/H$ is the $C^*$-subalgebra  of
$B(\H_\pi\otimes L^2(G))$ generated by the operators \[
\big\{(\pi\otimes\iota)\circ\delta(b)(1\otimes M|(f)):b\in B, f\in
C_0(G/H)\big\},\]
where $M|$ denotes the extension of $M$ to $C_0(G/H)\subset M(C_0(G))$. (The reasons
for calling this the \emph{reduced} crossed product are discussed in \cite[\S2]{EKR}.)
It  follows from
\cite[Lemma~11]{man} that 
\[ B\times_{\delta,r} 
(G/H)=\clsp\big\{(\pi\otimes\iota)\circ\delta(b)(1\otimes M|(f)):b\in B, f\in 
C_0(G/H)\big\}. \]

Suppose that $E$ is a compact subset of $G$ and $u\in A_c(G)$. Denote by $C_E(G)$ the set
of functions in
$C_c(G)$ with support in $E$.  Following Mansfield, we say that an operator on
$\H_\pi\otimes L^2(G)$ is $(u, E)$ if it is the norm limit of a sequence
$\{x_i\}$ of the form
\begin{equation}\label{Dnotation}
x_i=\sum_{j=1}^{n_i}
\big((\pi\otimes\iota)(\delta(\delta_u(b_{ij}))\big)\big(1\otimes
M(f_{ij})\big)\text{\ where\ } f_{ij}\in C_E(G).
\end{equation}
Then $\D$ denotes the
set of operators which are $(u, E)$ for some
$u\in A_c(G)$ and compact subset $E$ of $G$ \cite[\S 3]{man}.
The point of the definition is that one can pull 
$\pi\otimes\iota(\delta(\delta_u(b)))$
past  $(1\otimes M(f))$ and remain within $\D$  \cite[Lemma~9]{man} so that $\D$ is a
$*$-subalgebra of $B\times_\delta G$ \cite[Lemma~11]{man};  if $\delta$ is
nondegenerate, then  $\D$ is dense in $B\times_{\delta} G$.  

\section{The main theorem}\label{guts}

\begin{thm}\label{mainthm}
Suppose that $\delta$ is a nondegenerate reduced coaction of a locally compact group $G$ 
on a $C^*$-algebra $B$ and that $H$ is a closed subgroup of $G$.
Then the restriction of the dual action $\hat\delta: G\to\Aut(B\times_\delta G)$ to $H$ is
proper and saturated with respect to Mansfield's subalgebra $\D$,
and the generalised fixed-point algebra is the crossed product
$B\times_{\delta,r} (G/H)$ of $B$ by the homogeneous space $G/H$.
Thus, by \cite[Corollary~1.7]{rie-pr}, $\D$ completes to give a $((B\times_\delta
G)\times_{\hat\delta, r}H)$--$(B\times_{\delta,r} (G/H))$ imprimitivity bimodule.
\end{thm}

When the subgroup
$H$ is normal and amenable, all of this was proved by Mansfield in
\cite{man}, though he proved directly that $\D$ completes to give an imprimitivity
bimodule, and deduced from his arguments that the dual action of $H$ is proper. For
general $H$, it is known that $(B\times_\delta
G)\times_{\hat\delta, r}H$ and $B\times_{\delta,r} (G/H)$ are Morita equivalent
\cite[\S5]{EKR}; however, the argument in \cite{EKR} is indirect, and
\cite[Theorem~5.2]{EKR} generalises the weaker version of Mansfield's theorem proved in
\cite{Ng} rather than Mansfield's theorem itself. Theorem~\ref{mainthm} generalises
Mansfield's theorem, and gives a concrete bimodule which can be used to induce
representations.

The idea of approaching Mansfield's theorem through properness comes from
\cite{ahrw-symimp} and \cite{dpr}, where Theorem~\ref{mainthm} was proved for discrete
groups (\cite[Theorem~5.1]{dpr}). We know from recent work of Rieffel
\cite[Theorem~5.7]{rie-pr2} that $\hat\delta$ is proper and saturated 
with respect to the $*$-subalgebra
\begin{equation}
\label{rie-alg}
A_0:=(1\otimes M(C_c(G))(B\times_\delta G)(1\otimes M(C_c(G))
\end{equation} of $B\times_\delta G$ (see \cite[Remark~4.5]{ahrw-symimp}). It is 
not obvious to us, however, how to directly identify the 
generalised fixed-point algebra associated to $A_0$.  Our main observation is 
that if we use $\D$ instead of $A_0$, then we can use Mansfield's results to 
help identify the generalised fixed-point algebra with 
$B\times_{\delta,r}(G/H)$. The catch is that we cannot then use the results of 
\cite{rie-pr2}, and have to verify directly that $\hat\delta$ is proper with 
respect  to $\D$. 

Our first lemma says that $\D$ is smaller than Rieffel's $A_0$.

\begin{lemma}\label{rewrite}
Suppose $u\in A_c(G)$ and $E$ is a compact subset of $G$. Then there 
exist $f,g\in C_c(G)$
such that $x=(1\otimes M(f)) x (1\otimes M(g))$ for every $x\in \D$ 
which is $(u,E)$;
in particular, $\D=(1\otimes M(C_c(G)))\D(1\otimes M(C_c(G)))$.
\end{lemma}

\begin{proof}
Choose $g\in C_c(G)$ with $g=1$ on $E$. Then a glance at 
\eqref{Dnotation} shows that $x=x(1\otimes
M(g))$ whenever  $x$ is $(u, E)$. Let $\tilde u(t)=\overline{u(t^{-1})}$; then
by \cite[Lemma 11(ii)]{man}, there is a compact subset $F$ such that 
$x^*$ is $(\tilde u, F)$
whenever $x$ is $(u, E)$.  Now choose $f\in C_c(G)$ such
that $f=1$ on $F$, and we have $(1\otimes M(f))x=x$.
\end{proof}

To establish that $\hat\delta$ is proper, we verify the conditions of 
\cite[Definition~1.2]{rie-pr}:
for every $x,y\in\D$,
\begin{enumerate}
\item[1.]
both $s\mapsto x\hat\delta_s(y^*)$ and $s\mapsto
\Delta(s)^{-1/2}x\hat\delta_s(y^*)$ are integrable, and
\item[2.]  there exists $\langle x\,,\,y\rangle\in M(B\times_\delta 
G)$ such that for all $z\in\D$,
\begin{equation}\label{eq-invariant}
z\langle x\,,\,y\rangle\in\D \text{ \ and \ }z\langle x\,,\, y\rangle=\int_H 
z\hat\delta_s(x^*y)\, ds. \end{equation}
\end{enumerate}
Given Lemma~\ref{rewrite}, item 1 is easy: factor 
$x=x(1\otimes M(f))$ and $y=y(1\otimes M(g))$, so
that
\[
x\hat\delta_s(y^*)=x(1\otimes M(f\rt_s(\bar g))\hat\delta_s(y^*),
\]
and note that $s\mapsto f\rt_s(\bar g)$ has compact support.

So the crux in proving that $\hat\delta$ is proper is to find the 
multipliers $\langle
x\,,\,y\rangle$. To do this, we use an averaging technique which was 
developed by Olesen and Pedersen
\cite{OP1, OP2}, applied to coactions by Quigg \cite{quigg, QR}, and 
shown to be relevant to
properness by Rieffel \cite{rie-pr2}. Let $\alpha$ be an action of a
locally compact group $K$ on a $C^*$-algebra $A$. Denote by $\p$ the 
set of multipliers $a\in
M(A)^+$ for which there exists $\Psi(a)\in M(A)^+$ satisfying
\begin{equation}\label{charPsi}
\omega(\Psi(a)) =\int_K\omega(\alpha_s(a))\, ds \text{\ for all\
}\omega\in A^*
\end{equation}
(where we have implicitly extended $\omega$ to a strictly continuous 
functional on $M(A)$ and
$\alpha$ to a strictly continuous action on $M(A)\;$). The 
functionals $\omega\in A^*$ separate
points of
$M(A)$, so $\Psi(a)$ is uniquely determined by \eqref{charPsi}. By 
\cite[Corollary~3.6]{quigg}, the
linear span $\m$ of $\p$ is a $*$-subalgebra of $M(A)$ with 
$\m^+=\p$, and $\Psi$ extends uniquely
to a positive linear map $\Psi:\m\to M(A)$ such that $a\mapsto 
\Psi(fag)$ is norm-continuous on
$M(A)$ for fixed $f$ and $g$ in $\m$.

We now describe the averaging map $\Psi$ for
$(A,K,\alpha)=(B\times_\delta G,H,\hat\delta)$. For $f\in C_c(G)$, we 
define $\Phi(f)\in C_c(G/H)$
by
\[
\Phi(f)(tH)=\int_H f(ts)\,ds.
\]

\begin{lemma}\label{mlemma}
Rieffel's subalgebra 
$A_0$ of \eqref{rie-alg} is contained in $\m$, and in
particular $\D\subset\m$. For  $b\in 
B$, $u\in A_c(G)$ and $f\in
C_c(G)$ we have
\begin{equation}\label{pullbout}
\Psi\big((\pi\otimes\iota)(\delta(\delta_u(b)))(1\otimes M(f))\big)
=(\pi\otimes\iota)(\delta(\delta_u(b)))(1\otimes M(\Phi(f))).
\end{equation}
\end{lemma}

\begin{proof}
(These arguments are implicit in \cite[\S3]{quigg} and \cite[\S1]{QR}, see also 
\cite[\S4-5]{rie-pr2}.) We first claim that $f\in\m(C_0(G),\rt)$, with 
$\Psi(f)=\Phi(f)$.  To see this, it suffices to take
$f\in C_c(G)^+$, $\omega\in C_c(G)^{*+}$ and show that
\begin{equation}\label{psi=phi}
\omega(\Phi(f))=\int_H\omega(\rt_s(f))\, ds.
\end{equation}
But the Riesz Representation Theorem says that $\omega$ is given by a finite Borel
measure $\mu$, and then  \eqref{psi=phi} reduces to Tonelli's Theorem. 

Now \cite[Proposition~1.4]{QR} implies that $1\otimes M(f)$ belongs to
$\m(B\times_\delta G,H,\hat\delta)$; in fact, it belongs to the left 
ideal $\n$ of
\cite[Definition~3.4]{quigg}, and hence every element of the form 
$(1\otimes M(f))c(1\otimes
M(g))$ for $c\in B\times_\delta G$ belongs to $\m=\n^*\n$. 
Thus
$A_0\subset\m$.

To see \eqref{pullbout}, let $\kappa\in(B\times_\delta G)^*$ and define a
bounded functional $\omega$ on $C_0(G)$ by
 \[
 \omega(f):=\kappa \big(\pi\otimes\iota(\delta(\delta_u(b)))(1\otimes 
M(f))\big).  \]
Equation \eqref{psi=phi} extends to arbitrary bounded functionals by 
linearity, so 
\begin{align*}
 \kappa\big(\pi\otimes\iota(\delta(\delta_u(b))(1\otimes
M(\Phi(f)))\big)
&=\omega(\Phi(f))
=\int_H\omega(\rt_s(f))\, ds\\
&=\int_H\kappa\big(\pi\otimes\iota(\delta(\delta_u(b)))(1\otimes 
M(\rt_s(f)))\big)\, ds \\
&=\int_H\kappa\big(\hat\delta_s(\pi\otimes\iota(\delta(\delta_u(b)))(1\otimes
M(f)))\big)\, ds,
\end{align*}
and this implies \eqref{pullbout}.  \end{proof}

To finish the proof of Theorem~\ref{mainthm}, we need the subalgebra
$\D_H$ constructed in  \cite[\S 3]{man}; notice that Mansfield was careful not to use
amenability or normality in that section. Let $E$ be a compact subset
of $G$ and
$u\in A_c(G)$. An operator  on
$\H_\pi\otimes L^2(G)$ is $(u, E, H)$ if it is the norm limit of a
sequence $\{x_i\}$ of the form
\begin{equation}\label{dh}
x_i=\sum_{j=1}^{n_i}
\big((\pi\otimes\iota)(\delta(\delta_u(b_{ij}))\big)\big(1\otimes
M|(\Phi(f_{ij}))\big)\ \mbox{where  $f_{ij}\in C_E(G)$}.
\end{equation}
Then $\D_H$
is the set of operators on  $\H_\pi\otimes L^2(G)$ which are $(u, E, H)$ 
for some $u\in A_c(G)$ and compact subset $E$ of $G$.

\begin{proof}[The end of the proof of Theorem~\ref{mainthm}]
Given our earlier observations, to show that $\hat\delta$ is a proper action we 
need to define  multipliers  $\langle x\,,\, y\rangle$ of $\D$ satisfying \eqref{eq-invariant}. 
We set $\langle x\,,\, y\rangle:=\Psi(x^*y)$. 

We first claim that $\langle 
x\,,\, y\rangle\in\D_H$; this gives $z\langle x\,,\, y\rangle\in\D$ for 
$z\in\D$ because $\D_H$ multiplies $\D$ \cite[Lemma~11(v)]{man}.  Since $\D$ 
is a $*$-algebra we can write $x^*y$ as the norm limit of  \[  
x_i=\sum_{j=1}^{n_i} 
\big((\pi\otimes\iota)(\delta(\delta_u(b_{ij}))\big)\big(1\otimes M(f_{ij})\big)
 \text{\ where each $x_i$ is $(u,E)$.}
\] 
By Lemma~\ref{rewrite} there are  $f$, $g\in C_c(G)$  such that 
$x_i=(1\otimes M(f))x_i(1\otimes M(g))$ for all $i$, and then $x^*y=(1\otimes 
M(f))x^*y(1\otimes M(g))$ as well. The point is that since $1\otimes 
M(f)$ and $1\otimes 
M(g)$ belong to $\m$, the map $c\mapsto \Psi\big((1\otimes 
M(f))c(1\otimes M(g))\big)$  is norm-continuous on $B\times_\delta G$ by 
\cite[Corollary~3.6]{quigg}.  Thus we have  
\begin{align*} 
\langle x\,,\, y\rangle=\Psi(x^*y) 
&=\Psi( (1\otimes M(f))(\lim_i x_i)(1\otimes M(g))  )\\
&=\lim_i \Psi( (1\otimes M(f))(x_i)(1\otimes M(g))  )\\
&=\lim_i\Psi\Big(\sum_{j=1}^{n_i}
\big((\pi\otimes\iota)(\delta(\delta_u(b_{ij}))\big)\big(1\otimes
M(f_{ij})\big)\Big),
\end{align*}
and it follows from Lemma~\ref{mlemma} that
\begin{equation}
\langle x\,,\, y\rangle=\lim_i\sum_{j=1}^{n_i}
\big((\pi\otimes\iota)(\delta(\delta_u(b_{ij}))\big)\big(1\otimes
M(\Phi(f_{ij}))\big).\label{inDH}
\end{equation}
Thus $\langle x\,,\, y\rangle$  is $(u, E,H)$, and belongs to $\D_H$, as
claimed.

Next suppose $z\in\D$ and $\kappa\in (B\times_\delta G)^*$, and define a bounded
linear functional $\omega$ on $B\times_\delta G$ by $\omega(a)=\kappa(za)$.  Equation
\eqref{charPsi} gives
\begin{equation}\label{unique}
\kappa(z\langle x\,,\, y\rangle)=\omega(\Psi(x^*y))=\int_H
\omega(\hat\delta_s(x^*y))\, ds=\int_H\kappa(z\hat\delta_s(x^*y))\, ds,
\end{equation}
which equals $\kappa\big(\int_H z\hat\delta_s(x^*y))\,ds\big)$ because the 
integrand is integrable and bounded linear maps pull through integrals.  So 
$z\langle x\,,\, y\rangle=\int_H z\hat\delta_s(x^*y))\,ds$, and we have found the 
right multiplier. We have now proved that $\hat\delta$ is a proper action with respect to
$\D$.

A glance at \eqref{inDH} shows that each $\langle x\,,y\rangle$ belongs to 
$B\times_{\delta,r} (G/H)$; we  need to argue that  $\sp 
\{\langle x\,,\, y\rangle:x,y\in \D\}$ is dense in  
$B\times_{\delta,r} (G/H)$. Since $\delta$ is a nondegenerate coaction,  
$\{\delta_u(b):u\in A_c(G), b\in B\}$  is dense in $B$, and by linearity it 
suffices to show that each $((\pi\otimes\iota)(\delta(\delta_u(b))))(1\otimes 
M|(\Phi(f)))$ can be approximated by some $\langle 
x\,,\, y\rangle$.

Since $\D^2$ is dense in $B\times_\delta G$ there exist $x_i, y_i\in 
\D$ such that
\begin{equation}\label{cvgence}
x_i^*y_i\to (\pi\otimes\iota)(\delta(\delta_u(b)))(1\otimes M(f)).
\end{equation}
By applying Lemma~\ref{rewrite} to the right-hand side of (\ref{cvgence}), we may assume
there exist
$g_1,g_2\in C_c(G)$ such that 
\[
(1\otimes M(g_1))(x_i^*y_i)(1\otimes 
M(g_2))\to(\pi\otimes\iota)(\delta(\delta_u(b)))(1\otimes M(f)).
\]
Lemma~\ref{mlemma} and the continuity property of $\Psi$ 
\cite[Corollary~3.6]{quigg}  give 
\begin{align*}
(\pi\otimes\iota)(\delta(\delta_u(b)))(1\otimes M|(\Phi(f))
&=
\Psi((\pi\otimes\iota)(\delta(\delta_u(b)))(1\otimes M(f))\\
&=\Psi\big(((1\otimes M(g_1))\big(\lim_i x^*_iy_i\big)(1\otimes M(g_2))\big)\\
&=\lim_i\Psi\big(((1\otimes M(g_1))( x^*_iy_i)(1\otimes M(g_2))\big)\\
&=\lim_i\big\langle x_i(1\otimes M(g_1)^*)\,,\, y_i(1\otimes M(g_2))\big\rangle.
\end{align*}
Since $x_i(1\otimes M(g_1)^*)$ and  $y_i(1\otimes M(g_2))$ belong to $\D$ by
Lemma~\ref{rewrite}, this proves
that the generalised fixed-point algebra is all of
$B\times_{\delta,r} (G/H)$.

That $\hat\delta$ is also saturated with respect to $\D$ follows from the
proof of \cite[Lemma~4.1]{ahrw-symimp}: the proof goes through  if
$B\times_\delta G$ is replaced with any dense $*$-subalgebra.
The last assertion now follows directly from \cite[Corollary~1.7]{rie-pr}. 
\end{proof}

\begin{remark}\label{remark-rie}
It follows from our result   that the generalised fixed-point algebra 
associated to Rieffel's $*$-subalgebra $A_0:=(1\otimes M(C_c(G))(B\times_\delta G)(1\otimes M(C_c(G))$ 
coincides with $B\times_{\delta, r}(G/H)$.  To see this, let $a,b\in A_0$. The 
calculation~\eqref{unique} shows that Rieffel's inner product 
$\langle a\,,\, b\rangle_D$ coincides with $\Psi(a^*b)$.  Now choose 
$f,g\in C_c(G)$ such that $a^*b=(1\otimes M(f))a^*b(1\otimes M(g))$, 
and sequences $\{ x_i\}$, $\{y_i\}$ in $\D$ such that $x_i\to a$ and 
$y_i\to b$.   The continuity property of $\Psi$ 
described in \cite[Corollary~3.6]{quigg} gives
\begin{align*}
\Psi(a^*b)&=\Psi\big(\lim_i (1\otimes M(f))x_i^*y_i(1\otimes M(g)) \big)\\
&=\lim_i\Psi\big( (1\otimes M(f))x_i^*y_i(1\otimes M(g))\big)\\
&=\lim_i\big\langle x_i(1\otimes M(f))\,,\, y_i(1\otimes M(g))\big\rangle,
\end{align*}
so $\Psi(a^*b)$ is the limit of elements in $B\times_{\delta,r}(G/H)$ and is 
itself in   $B\times_{\delta,r}(G/H)$.   Since $\{\langle x\,,\, 
y\rangle:x,y\in\D\}$ spans a dense subspace of $B\times_{\delta,r}(G/H)$,  so does 
$\{\langle  a\,,\, b\rangle_D:a,b\in A_0\}$. 
\end{remark}

If $\nu$ is a nondegenerate representation of $B\times_{\delta,r}(G/H)$ on 
$\H$, we can induce it to a representation $\overline{\D}\dashind\nu$ of 
$(B\times_{\delta} G)\times_{\hat\delta, r}H$ on 
$\overline{\D}\otimes_{B\times_{\delta,r}(G/H)}\H$. Composing 
$\overline{\D}\dashind\nu$ with the canonical map of $B\times_\delta G$ into 
$M((B\times_\delta G)\times_{\hat\delta}H)$ gives a representation 
$\Ind_{G/H}^G\nu$ of $B\times_\delta G$ which we can reasonably call the 
\emph{representation of $B\times_\delta G$ induced from $\nu$}.  There is also a 
unitary representation $U$ of $H$ such that $(\Ind_{G/H}^G\nu, U)$ is a 
covariant representation of $(B\times_\delta G, H,\hat\delta)$, and the 
representation $(\Ind_{G/H}^G\nu)\times  U$ of $(B\times_\delta 
G)\times_{\hat\delta} H$ factors through the representation $\overline{\D}\dashind\nu$ of
the reduced crossed  product.   If $H$ is amenable then  $(B\times_\delta 
G)\times_{\hat\delta} H=(B\times_\delta 
G)\times_{\hat\delta,r} H$, and it follows from the Rieffel correspondence that  
every representation  $\mu\times U$ of $(B\times_\delta G)\times_{\hat\delta}
H$ is equivalent to one of the form  $\overline{\D}\dashind \nu$. Thus we 
have the following imprimitivity theorem:

\begin{cor} 
\label{cor-amenable} Suppose  $\delta$ is a nondegenerate reduced coaction of a locally 
compact group $G$ on a $C^*$-algebra $B$,  $H$ is a closed amenable subgroup of 
$G$, and $\mu$ is a nondegenerate representation of $B\times_\delta G$ on $\H$.  
Then  $\mu$  is unitarily equivalent to a representation induced 
 from a representation $\nu$ of $B\times_{\delta,r}(G/H)$ if 
and only if there is a unitary representation $U$
of $H$ on $\H$ such that $(\mu, U)$ is a covariant representation of
$(B\times_\delta G, H,\hat\delta)$.
  \end{cor}

\section{Normal subgroups} \label{normal}

Throughout this section, $N$ is a closed normal subgroup of $G$, so that
$G/N$ is itself a locally compact group. We will show how to
restrict the coaction $\delta$ to a coaction $\delta|$ of $G/N$, and that the
crossed product by $\delta|$ is   then   isomorphic to $B\times_{\delta,
r}(G/N)$. We therefore obtain an imprimitivity theorem relating $(B\times_\delta
G)\times_{\hat\delta,r}N$ to $B\times_{\delta|}(G/N)$. From this we can deduce a
theorem of Quigg and Kaliszewski \cite[Theorem~3.3]{KQ} which extends Mansfield
imprimitivity to full coactions and normal nonamenable subgroups. Our main
innovation here is the concept of restriction for reduced coactions, which is
 based on work of Quigg \cite{quigg-fr}.

Let $\delta:B\to M(B\otimes C^*_r(G))$ be a nondegenerate coaction.
The quotient map  $q:G\to G/N$ induces a homomorphism $q:C^*(G)\to
C^*(G/N)$. The composition $\lambda^{G/N}\circ q$ is the induced representation
$\Ind_N^G (1^N)$, which satisfies $\ker \Ind_N^G (1^N)\supset
\ker\Ind(\lambda^N)=\ker\lambda^G$ if and only if $\ker 1^N\supset
\ker\lambda^N$, and hence if and only if $N$ is amenable. If $N$ is amenable,
therefore, $q$ descends to a homomorphism $q^r:C^*_r(G)\to C^*_r(G/N)$, and the
restriction is usually defined to be $\delta|:=(\iota\otimes q^r)\circ\delta$.
For this to work, however, we do not need that $\lambda^{G/N}\circ q$ factors
through $C^*_r(G)$, only that $\iota\otimes (\lambda^{G/N}\circ q)$ factors
through $B\otimes C^*_r(G)$,  and this happens more frequently than one 
might expect: for example, it follows from the absorbing property of the regular 
representation that this happens for $B=C^*_r(G)$. So there remains the 
possibility of defining restrictions even when $N$ is not amenable.
The usual way to get round this problem of defining restrictions is to consider
only full coactions $\epsilon:B\to M(B\otimes C^*(G))$, because then
$\epsilon|:=(\iota\otimes q)\circ \epsilon:B\to M(B\otimes C^*(G/N))$ always
makes sense. 

A full coaction is \emph{Quigg-normal} if the canonical map $j_B$ of $B$
into $M(B\times_\epsilon G)$ is injective 
\cite[Definition~2.1]{quigg-fr}.  Quigg showed in \cite[Theorem~4.7]{quigg-fr} 
that every nondegenerate reduced coaction $\delta$ of $G$ on $B$ admits a unique 
\emph{Quiggification}: a  nondegenerate full Quigg-normal  coaction 
$\delta^q:B\to M(B\otimes C^*(G))$ with reduction 
$(\delta^q)^r:=(\iota\otimes\lambda_G)\circ\delta^q$ equal to $\delta$. If $N$ 
is amenable, the equation $\lambda^{G/N}\circ q=q^r\circ \lambda^G$ implies that 
$\delta|=(\delta^q)^r|=(\delta^q|)^r$. Thus we feel free to define:

\begin{definition}  \label{def-restrict}
 Let $\delta:B\to M(B\otimes C^*_r(G))$ be a  coaction of $G$ on
$B$. Then  the \emph{restriction}
of $\delta$ to a coaction of $G/N$ is
$\delta|:=(\delta^q|)^r: B\to M(B\otimes C^*_r(G/N))$.
\end{definition}

The crucial points about this definition are, first, that it makes sense for
any normal subgroup, amenable or not, and, second, that the crossed product by
$\delta|$ coincides with  $B\times_{\delta,r}(G/N)$. To see this second
point, we consider the canonical covariant representation $(j_B,j_G)$ of
$(B,G,\delta^q)$ in $M(B\times_{\delta^q}G)$, and the restriction $j_G|$ of
$j_G$ to $C_0(G/N)\subset M(C_0(G))$. Then $(j_B,j_G|)$ is a covariant
representation of $(B,G/N,\delta^q|)$ in $M(B\times_{\delta^q}G)$, and because
$\delta^q$ is Quigg-normal, \cite[Lemma~3.2]{KQ} implies that
$j_B\times j_G|$ is an injective homomorphism of $B\times_{\delta^q|}(G/N)$ into
$M(B\times_{\delta^q}G)$; its range is spanned by elements of the form
$j_B(b)j_G(f)$ for $f\in C_0(G/N)$. Since
$B\times_{\delta^q}G$ is canonically isomorphic to
$B\times_{(\delta^q)^r}G=B\times_\delta G$ \cite[Theorem~4.1]{R}, and this
isomorphism carries $\clsp\{j_B(b)j_G(f):f\in C_0(G/N)\}$ onto
$B\times_{\delta,r}(G/N)$, we have:

\begin{prop}\label{prop-allOK}
  Let $\delta$
be  a nondegenerate reduced coaction of $G$ on $B$, let $N$ be a closed normal
subgroup of $G$, and let $\delta|$ be the restriction to a reduced coaction of
Definition~\ref{def-restrict}. Then $B\times_{\delta|}(G/N)$  is isomorphic to
the crossed product $B\times_{\delta,r}(G/N)$ of $B$ by the homogeneous space
$G/N$. \end{prop}

\begin{cor}\label{cor-normal}
Let $\delta$ be  a nondegenerate reduced coaction of $G$ on $B$. Then
Mansfield's subalgebra $\D$ of $B\times_\delta G$ completes to give a Morita equivalence
between
$(B\times_\delta G)\times_{\hat\delta, r} N$ and
$B\times_{\delta|}(G/N)$. \end{cor}

Using the trick of passing from full to reduced crossed products 
we obtain a new proof of the 
Kaliszewski-Quigg version of Mansfield's imprimitivity theorem using 
Corollary~\ref{cor-normal}.

\begin{cor}[{\cite[Corollary~3.4]{KQ}}]\label{cor-kq}
Suppose $\epsilon:B\to
M(B\otimes C^*(G))$ is a nondegenerate full Quigg-normal coaction and $N$ is a
closed normal subgroup of $G$. Then $\D$ completes to give a Morita equivalence
between $(B\times_\epsilon G)\times_{\hat\epsilon, r}N$ and
$B\times_{\epsilon|}(G/N)$.
 \end{cor}

 \begin{proof}
Let $I=\ker((\iota\otimes\lambda)\circ\epsilon)$, and consider the reduction
$\epsilon^r:B/I\to M((B/I)\otimes C^*_r(G))$.  It follows from Corollary~\ref{cor-normal}
that
$\D$ completes to give a Morita equivalence between $((B/I)\times_{\epsilon^r}
G)\times_{\widehat{\epsilon^r},r}N$ and $(B/I)\times_{\epsilon^r|}
(G/N)$. By \cite[Theorem~4.1]{R} there is an isomorphism $\Q$ of
$(B/I)\times_{\epsilon^r}G$ onto $B\times_\epsilon G$ such that
\[
\Q( j_{B/I}^r(q(b))j_G(f))=j_B(b)j_G(f).
\]
Thus $\Q$ preserves the dual actions, and induces an isomorphism of
$((B/I)\times_{\epsilon^r}G)\times_{\widehat{\epsilon^r},r}N$ onto
$(B\times_\epsilon G)\times_{\hat\epsilon, r}N$.
In addition, $\Q$ carries
  \[
  j_{B/I}^r(q(b))j_G|(f)\text{\ to\ }j_B(b)j_G|(f),
  \]
  and hence is an isomorphism of $(B/I)\times_{\epsilon^r|} (G/N)$ onto
$B\times_{\epsilon|} (G/N)$.
  \end{proof}

\begin{remark}\label{kq}
Because we have made sense of restriction for reduced coactions, we 
can also prove Corollary~\ref{cor-normal} by applying Corollary~\ref{cor-kq} to 
the Quiggification $\delta^q$, which is  always Quigg-normal. 
\end{remark}


\begin{thebibliography}{20}

\bibitem{dpr}
K. Deicke, D. Pask, and I. Raeburn,
\emph{Coverings of directed graphs and crossed products of 
$C^*$-algebras by coactions of homogeneous
spaces}, preprint (arXiv.math.OA/0201033), 2001.

\bibitem{EKR} S. Echterhoff, S. Kaliszewski, and I. Raeburn,
\emph{Crossed products by dual coactions of groups and homogeneous spaces},
J. Operator Theory \textbf{39} (1998), 151--176.

\bibitem{EQ} S. Echterhoff and J. Quigg, \emph{Full duality for coactions of discrete
groups}, Math. Scand., to appear.

\bibitem{green} P. Green, \emph{The local structure of twisted covariance
algebras}, Acta Math. \textbf{140} (1978), 191--250.

\bibitem{ahrw-symimp} A. an Huef, I. Raeburn, and D.P. Williams,
\emph{A symmetric imprimitivity theorem for commuting proper actions}, preprint
(arXiv.math.OA/0202046), 2002.

\bibitem{KQ} S. Kaliszewski and J. Quigg, \emph{Imprimitivity for
$C^*$-coactions of non-amenable groups}, Math. Proc. Camb. Phil. Soc.
\textbf{123} (1998), 101--118.

\bibitem{lprs} M.B. Landstad, J. Phillips, I. Raeburn and C.E. Sutherland,
\emph{Representations of crossed products by coactions and principal bundles}, Trans.
Amer. Math. Soc. \textbf{299} (1987), 747--784.

\bibitem{man} K. Mansfield, \emph{Induced
representations of crossed products by coactions}, J. Funct. Anal.
\textbf{97}   (1991), 112--161.

\bibitem{Ng} C.K. Ng, \emph{A remark on Mansfield's imprimitivity theorem}, 
Proc. Amer. Math. Soc. \textbf{126} (1998), 3767--3768.

  \bibitem{OP1} D. Olesen and G.K. Pedersen, \emph{Applications of the Connes
spectrum to $C^*$-dynamical systems}, J. Funct. Anal. \textbf{30} (1978),
179--197.

\bibitem{OP2} D. Olesen and G.K. Pedersen, \emph{Applications
of the Connes spectrum to $C^*$-dynamical systems \textrm{II}},
J. Funct. Anal. \textbf{36} (1980), 18--32.


\bibitem{quigg} J.C. Quigg,
\emph{Landstad duality for $C^*$-coactions},  Math. Scand. \textbf{71} (1992),
277--294.

\bibitem{quigg-fr}  J.C. Quigg,
\emph{Full and reduced coactions}, Math. Proc. Camb. Phil. Soc. \textbf{116}
(1994), 435--450.

  \bibitem{QR} J.C. Quigg and I. Raeburn, \emph{Induced
$C^*$-algebras and Landstad duality for twisted coactions}, Trans. Amer. Math.
Soc. \textbf{347} (1995), 2885--2915.

\bibitem{R} I. Raeburn, \emph{On crossed products by coactions and 
their representation theory},
Proc. London Math. Soc. \textbf{64} (1992), 625--652.

\bibitem{rie-pr} M.A. Rieffel, \emph{Proper actions of groups on
     $C^*$-algebras}, Mappings of Operator Algebras, Progr. Math., 
vol.~84, Birkhauser, Boston,
1988, pp. 141--182.


\bibitem{rie-pr2} M.A. Rieffel, \emph{Integrable and proper actions on
     $C^*$-algebras, and square integrable representations of groups},
   preprint (arXiv.math.OA/9809098), 1999.


\end{thebibliography}
\end{document}